\documentclass{article}
\usepackage{amsmath, amsthm, amsfonts}

\theoremstyle{definition}

\theoremstyle{remark}

\input{tcilatex}
\begin{document}

\title{An elliptic partial differential equations system and its application}
\author{Dragos-Patru Covei \\
{\small The Bucharest University of Economic Studies, Bucharest, Romania. }\\
{\small E-mail address: dragos.covei@csie.ase.ro} \and Traian A. Pirvu \\
{\small McMaster University, Hamilton, Canada. }\\
{\small E-mail address: tpirvu@math.mcmaster.ca}}
\maketitle

\abstract{
This paper deals with the following elliptic system of equations 
\begin{equation}
\left\{ 
\begin{array}{c}
-\frac{k_{1}}{2}\Delta z_{1}(x)+\frac{|\nabla z_{1}(x)|^{2}}{2}%
=f_{1}(x)-(\lambda _{1}+a_{1})z_{1}(x)+a_{1}z_{2}(x), \\ 
-\frac{k_{2}}{2}\Delta z_{2}(x)+\frac{|\nabla z_{2}(x)|^{2}}{2}%
=f_{2}(x)-(\lambda _{2}+a_{2})z_{2}(x)+a_{2}z_{1}(x),%
\end{array}%
\right. x\in \mathbb{R}^{N}\text{,}
\end{equation}%
where $\lambda _{i}>0$ ($i=1,2$) are some real constants suitable chosen and 
$k_{i}>0$, $a_{i}>0$ ($i=1,2$) are some real arbitrary constants, and $f_{i}$
are some continuous functions. The solution method is based on the sub- and
super-solutions approach. An application to a stochastic control problem is
presented. This system seemed not considered before.}

\section{Introduction}

In this work we study the existence of positive solutions for the following
partial differential equations (PDE) system

\begin{equation}
\left\{ 
\begin{array}{c}
-\frac{k_{1}}{2}\Delta z_{1}(x)+\frac{|\nabla z_{1}(x)|^{2}}{2}%
=f_{1}(x)-(\lambda _{1}+a_{1})z_{1}(x)+a_{1}z_{2}(x), \\ 
-\frac{k_{2}}{2}\Delta z_{2}(x)+\frac{|\nabla z_{2}(x)|^{2}}{2}%
=f_{2}(x)-(\lambda _{2}+a_{2})z_{2}(x)+a_{2}z_{1}(x),%
\end{array}%
\right. x\in \mathbb{R}^{N}\text{.}  \label{tpes}
\end{equation}%
Here $N\geq 1$ is the space dimension, $\left\vert \circ \right\vert $ is
the Euclidean norm of $\mathbb{R}^{N}$, $f_{i}:\mathbb{R}^{N}\rightarrow %
\left[ 0,\infty \right) $ ($i=1,2$) are continuous convex functions
satisfying%
\begin{equation}
\text{there exists }M_{i}>0\text{ such that }f_{i}\left( x\right) \leq
M_{i}\left( \left\vert x\right\vert ^{2}+1\right) ,  \label{acond}
\end{equation}%
$\lambda _{i}>0$ ($i=1,2$) are some real constants suitable chosen and $%
k_{i}>0$, $a_{i}>0$ ($i=1,2$) are some real arbitrary constants.

This system has received much attention in the last decades since it is
related with several models that arise in different mathematical models of
natural phenomena; for more on this see the papers of Akella and Kumar \cite%
{AK}, Alvarez \cite{A}, Bensoussan, Sethi, Vickson and Derzko \cite{BS}, the
authors \cite{CP}, Ghosh, Arapostathis and Marcus \cite{GAS} and Lasry and
Lions \cite{LS}.

The system studied in this paper appears naturally in characterizing value
function of a stochastic control problem with regime switching. Regime
switching is a phenomena present in many real world problems. Let us point
to some relevant works where regime switching is present. In finance \cite%
{PZ} studies the effect of regime switching (high versus low discount rates)
to a consumption and investment problem. In production management \cite{CLP}
studies the cost minimization problem of a company within an economy
characterized by two regimes. In civil engineering \cite{DMDK} studies the
optimal stochastic control problem for home energy systems with regime
switching; the two regimes are the peak and off peak energy demand.

The principal device in studying this system comes from the recent work of 
\cite{DPB}, where the author obtained non-positive radial solutions for the
system (\ref{tpes}) and where we postulate an open problem regarding the
existence of positive solution for this system. Another goal of this paper
is to improve the model given in \cite{BS}, \cite{DPB}, \cite{GAS} and to
give a verification result, i.e., show that the solution of the system
yields the optimal control.

Furthermore, there seems to be no previous mathematical results about the
existence of positive solutions for the semilinear system (\ref{tpes}). This
should not surprise us since there are some difficulties in analyzing this
class of systems in $\mathbb{R}^{N}$ ($N\geq 1$), which will be revealed in
the following sections organized as follows. In Section \ref{mr}, we give
our main theorem regarding the existence of positive solution for the
problem (\ref{tpes}) and its proof. Section \ref{ppp} contains the context
and the diffusion model from where such system appear. Section \ref{avr},
presents a verification result. In Section \ref{sc} we obtain a closed form
solution for our system in a special case.

\section{Main Result \label{mr}}

The main result of our paper, a basic existence theorem for (\ref{tpes}), is
presented.

\begin{theorem}
\label{msg}For all $\lambda _{1},\lambda _{2}\in \left( 0,\infty \right) $
the system of equations (\ref{tpes}) has a unique positive classical convex
solution with quadratic growth, i.e., 
\begin{equation}
z_{i}\left( x\right) \leq K_{i}(1+|x|^{2})\text{, for some }K_{i}>0,\quad
i=1,2,  \label{ineq}
\end{equation}%
and, such that%
\begin{equation}
\left\vert \nabla {z}_{i}(x)\right\vert \leq \overline{C}_{i}(1+\left\vert
x\right\vert ),\text{ for }x\in \mathbb{R}^{N}\text{and for some positive
constant }\overline{C}_{i}.  \label{sass}
\end{equation}
\end{theorem}

We give a detailed proof of Theorem \ref{msg}, which is based on the
following two results.

\begin{lemma}
\label{ech}The system of partial differential equations with gradient term (%
\ref{tpes}) is equivalent to the semilinear elliptic system%
\begin{equation}
\left\{ 
\begin{array}{l}
\Delta u=u\left( x\right) [\frac{2}{k_{1}^{2}}\left( f_{1}\left( x\right)
+\left( \lambda _{1}+a_{1}\right) k_{1}\ln u-a_{1}k_{2}\ln v\right) ], \\ 
\Delta v=v\left( x\right) [\frac{2}{k_{2}^{2}}\left( f_{2}\left( x\right)
+\left( \lambda _{2}+a_{2}\right) k_{2}\ln v-a_{2}k_{1}\ln u\right) ],%
\end{array}%
\right. \text{ }x\in \mathbb{R}^{N}.  \label{tpe}
\end{equation}
\end{lemma}

\paragraph{\textbf{{Proof}.}}

\bigskip The change of variable 
\begin{equation*}
z_{1}\left( x\right) =k_{1}w_{1}\left( x\right) \text{ and }z_{2}\left(
x\right) =k_{2}w_{2}\left( x\right) ,
\end{equation*}%
transform the system (\ref{tpes}) into%
\begin{equation}
\left\{ 
\begin{array}{c}
-\frac{k_{1}^{2}}{2}\Delta w_{1}+\frac{k_{1}^{2}\left\vert \nabla
w_{1}\right\vert ^{2}}{2}=f_{1}\left( x\right) -\left( \lambda
_{1}+a_{1}\right) k_{1}w_{1}+a_{1}k_{2}w_{2}, \\ 
-\frac{k_{2}^{2}}{2}\Delta w_{2}+\frac{k_{2}^{2}\left\vert \nabla
w_{2}\right\vert ^{2}}{2}=f_{2}\left( x\right) -\left( \lambda
_{2}+a_{2}\right) k_{1}w_{2}+a_{2}k_{1}w_{1},%
\end{array}%
\right.  \label{13}
\end{equation}%
or, equivalently%
\begin{equation}
\left\{ 
\begin{array}{c}
-\Delta w_{1}+\left\vert \nabla w_{1}\right\vert ^{2}=\frac{2}{k_{1}^{2}}%
\left[ f_{1}\left( x\right) -\left( \lambda _{1}+a_{1}\right)
k_{1}w_{1}+a_{1}k_{2}w_{2}\right] , \\ 
-\Delta w_{2}+\left\vert \nabla w_{2}\right\vert ^{2}=\frac{2}{k_{2}^{2}}%
\left[ f_{2}\left( x\right) -\left( \lambda _{2}+a_{2}\right)
k_{2}w_{2}+a_{2}k_{1}w_{1}\right] .%
\end{array}%
\right.  \label{sis2}
\end{equation}%
The change of variable 
\begin{equation*}
u\left( x\right) =e^{-w_{1}\left( x\right) }\text{ and }v\left( x\right)
=e^{-w_{2}\left( x\right) },
\end{equation*}%
transform the system (\ref{sis2}) into%
\begin{equation}
\left\{ 
\begin{array}{l}
\Delta u=u[\frac{2}{k_{1}^{2}}\left( f_{1}\left( x\right) +\left( \lambda
_{1}+a_{1}\right) k_{1}\ln u-a_{1}k_{2}\ln v\right) ], \\ 
\Delta v=v[\frac{2}{k_{2}^{2}}\left( f_{2}\left( x\right) +\left( \lambda
_{2}+a_{2}\right) k_{2}\ln v-a_{2}k_{1}\ln u\right) ],%
\end{array}%
\right.  \label{sisf}
\end{equation}%
since%
\begin{equation}
\begin{array}{cc}
\Delta u\left( x\right) =e^{-w_{1}\left( x\right) }(-\Delta w_{1}\left(
x\right) +\left\vert \nabla w_{1}\left( x\right) \right\vert ^{2}), &  \\ 
\Delta v\left( x\right) =e^{-w_{2}\left( x\right) }(-\Delta w_{2}\left(
x\right) +\left\vert \nabla w_{2}\left( x\right) \right\vert ^{2}). & 
\end{array}
\label{15}
\end{equation}

The existence of a solution $\left( u\left( x\right) ,v\left( x\right)
\right) \in C^{2}\left( \mathbb{R}^{N}\right) \times C^{2}\left( \mathbb{R}%
^{N}\right) $ for the problem (\ref{tpe}), such that $0<u\left( x\right)
\leq 1$ and $0<v\left( x\right) \leq 1$, for all $x\in \mathbb{R}^{N}$, is
proved in the following:

\begin{theorem}
\label{mss}If there exist functions $\underline{u}$, $\underline{v}$, $%
\overline{u}$, $\overline{v}:\mathbb{R}^{N}\rightarrow (0,1]$ of class $%
C^{2}\left( \mathbb{R}^{N}\right) $ such that 
\begin{equation}
\left\{ 
\begin{array}{l}
-\Delta \underline{u}\left( x\right) +\underline{u}\left( x\right) [\frac{2}{%
k_{1}^{2}}\left( f_{1}\left( x\right) +\left( \lambda _{1}+a_{1}\right)
k_{1}\ln \underline{u}\left( x\right) \right) ]\leq 2a_{1}\frac{k_{2}}{%
k_{1}^{2}}\underline{u}\left( x\right) \ln \underline{v}\left( x\right) , \\ 
-\Delta \underline{v}\left( x\right) +\underline{v}\left( x\right) [\frac{2}{%
k_{2}^{2}}\left( f_{2}\left( x\right) +\left( \lambda _{2}+a_{2}\right)
k_{2}\ln \underline{v}\left( x\right) \right) ]\leq 2a_{2}\frac{k_{1}}{%
k_{2}^{2}}\underline{v}\left( x\right) \ln \underline{u}\left( x\right) , \\ 
-\Delta \overline{u}\left( x\right) +\overline{u}\left( x\right) [\frac{2}{%
k_{1}^{2}}\left( f_{1}\left( x\right) +\left( \lambda _{1}+a_{1}\right)
k_{1}\ln \overline{u}\left( x\right) \right) ]\geq 2a_{1}\frac{k_{2}}{%
k_{1}^{2}}\overline{u}\left( x\right) \ln \overline{v}\left( x\right) , \\ 
-\Delta \overline{v}\left( x\right) +\overline{v}\left( x\right) [\frac{2}{%
k_{2}^{2}}\left( f_{2}\left( x\right) +\left( \lambda _{2}+a_{2}\right)
k_{2}\ln \overline{v}\left( x\right) \right) ]\geq 2a_{2}\frac{k_{1}}{%
k_{2}^{2}}\overline{v}\left( x\right) \ln \overline{u}\left( x\right) , \\ 
\underline{u}\left( x\right) \leq \overline{u}\left( x\right) \text{, \ }%
\underline{v}\left( x\right) \leq \overline{v}\left( x\right) ,%
\end{array}%
\right.  \label{tple}
\end{equation}%
in the entire Euclidean space $\mathbb{R}^{N}$, then system (\ref{tpe})
possesses an entire solution $\left( u,v\right) \in C^{2}\left( \mathbb{R}%
^{N}\right) \times C^{2}\left( \mathbb{R}^{N}\right) $ with $\underline{u}%
\left( x\right) \leq u\left( x\right) \leq \overline{u}\left( x\right) $ in $%
\mathbb{R}^{N}$ and $\underline{v}\left( x\right) \leq v\left( x\right) \leq 
\overline{v}\left( x\right) $ in $\mathbb{R}^{N}$.
\end{theorem}

Let us point out that the functions $\left( \underline{u},\underline{v}%
\right) $ (resp. $\left( \overline{u},\overline{v}\right) $) are called
sub-solution (resp. super-solution) for the system (\ref{tpe}).

\subparagraph{\textbf{{Proof}.}}

In the following we construct the functions $\left( \underline{u},\underline{%
v}\right) ,$ $\left( \overline{u},\overline{v}\right) $ which satisfies the
inequalities (\ref{tpe}) in $\mathbb{R}^{N}$. We proceed as in Bensoussan,
Sethi, Vickson and Derzko \cite{BS}, for the scalar case. More exactly, we
observe that there exist 
\begin{equation}
\left( \underline{u}\left( x\right) ,\underline{v}\left( x\right) \right)
=\left( e^{B_{1}\left\vert x\right\vert ^{2}+D_{1}},e^{B_{2}\left\vert
x\right\vert ^{2}+D_{2}}\right) ,\text{ with }B_{1},B_{2},D_{1},D_{2}\in
\left( -\infty ,0\right) ,  \label{subsuper}
\end{equation}%
such that for all $\lambda _{1}>0$ and \ $\lambda _{2}>0$ the following hold%
\begin{equation}
\left\{ 
\begin{array}{l}
-\Delta \underline{u}\left( x\right) +\underline{u}\left( x\right) [\frac{2}{%
k_{1}^{2}}\left( f_{1}\left( x\right) +\left( \lambda _{1}+a_{1}\right)
k_{1}\ln \underline{u}\left( x\right) \right) ]\leq 2a_{1}\frac{k_{2}}{%
k_{1}^{2}}\underline{u}\left( x\right) \ln \underline{v}\left( x\right) , \\ 
-\Delta \underline{v}\left( x\right) +\underline{v}\left( x\right) [\frac{2}{%
k_{2}^{2}}\left( f_{2}\left( x\right) +\left( \lambda _{2}+a_{2}\right)
k_{2}\ln \underline{v}\left( x\right) \right) ]\leq 2a_{2}\frac{k_{1}}{%
k_{2}^{2}}\underline{v}\left( x\right) \ln \underline{u}\left( x\right) ,%
\end{array}%
\right.  \label{lam}
\end{equation}%
i.e. $\left( \underline{u}\left( x\right) ,\underline{v}\left( x\right)
\right) $ is a sub-solution for the problem (\ref{tpe}). Indeed, we find $%
B_{1}$, $B_{2}$, $D_{1}$, $D_{2}\in \left( -\infty ,0\right) $ such that

\begin{equation*}
\left\{ 
\begin{array}{l}
-\allowbreak 2B_{1}\left( 2\left\vert x\right\vert ^{2}B_{1}-1\right)
-2B_{1}\left( N-1\right) +\frac{2}{k_{1}^{2}}\left[ M_{1}\left( \left\vert
x\right\vert ^{2}+1\right) +\left( \lambda _{1}+a_{1}\right) k_{1}\left(
B_{1}\left\vert x\right\vert ^{2}+D_{1}\right) \right]  \\ 
=\frac{2a_{1}k_{2}}{k_{1}^{2}}\left( B_{2}\left\vert x\right\vert
^{2}+D_{2}\right) , \\ 
-\allowbreak 2B_{2}\left( 2\left\vert x\right\vert ^{2}B_{2}-1\right)
-2B_{2}\left( N-1\right) +\frac{2}{k_{2}^{2}}\left[ M_{2}\left( \left\vert
x\right\vert ^{2}+1\right) +\left( \lambda _{2}+a_{2}\right) k_{2}\left(
B_{2}\left\vert x\right\vert ^{2}+D_{2}\right) \right]  \\ 
=\frac{2a_{2}k_{1}}{k_{2}^{2}}\left( B_{1}\left\vert x\right\vert
^{2}+D_{1}\right) ,%
\end{array}%
\right. 
\end{equation*}%
or, after rearranging the terms%
\begin{equation*}
\left\{ 
\begin{array}{l}
\left\vert x\right\vert ^{2}\left[ -4B_{1}^{2}+\frac{2M_{1}}{k_{1}^{2}}+%
\frac{2}{k_{1}^{2}}\left( \lambda _{1}+a_{1}\right) k_{1}B_{1}-2a_{1}\frac{%
k_{2}}{k_{1}^{2}}B_{2}\right] -2B_{1}N \\ 
+\frac{2M_{1}}{k_{1}^{2}}+\frac{2}{k_{1}^{2}}\left( \lambda
_{1}+a_{1}\right) k_{1}D_{1}-\frac{2a_{1}k_{2}D_{2}}{k_{1}^{2}}=0, \\ 
\left\vert x\right\vert ^{2}\left[ -4B_{2}^{2}+\frac{2M_{2}}{k_{2}^{2}}+%
\frac{2}{k_{2}^{2}}\left( \lambda _{2}+a_{2}\right) k_{2}B_{2}-2a_{2}\frac{%
k_{1}}{k_{2}^{2}}B_{1}\right] -2B_{2}N \\ 
+\frac{2M_{2}}{k_{2}^{2}}+\frac{2}{k_{2}^{2}}\left( \lambda
_{2}+a_{2}\right) k_{2}D_{2}-\frac{2a_{2}k_{1}D_{1}}{k_{2}^{2}}=0.%
\end{array}%
\right. 
\end{equation*}%
Now, we consider the system of equations%
\begin{equation}
\left\{ 
\begin{array}{l}
-4B_{1}^{2}+\frac{2M_{1}}{k_{1}^{2}}+\frac{2}{k_{1}^{2}}\left( \lambda
_{1}+a_{1}\right) k_{1}B_{1}-2a_{1}\frac{k_{2}}{k_{1}^{2}}B_{2}=0 \\ 
-2B_{1}N+\frac{2M_{1}}{k_{1}^{2}}+\frac{2}{k_{1}^{2}}\left( \lambda
_{1}+a_{1}\right) k_{1}D_{1}-2a_{1}\frac{k_{2}}{k_{1}^{2}}D_{2}=0 \\ 
-4B_{2}^{2}+\frac{2M_{2}}{k_{2}^{2}}+\frac{2}{k_{2}^{2}}\left( \lambda
_{2}+a_{2}\right) k_{2}B_{2}-2a_{2}\frac{k_{1}}{k_{2}^{2}}B_{1}=0 \\ 
-2B_{2}N+\frac{2M_{2}}{k_{2}^{2}}+\frac{2}{k_{2}^{2}}\left( \lambda
_{2}+a_{2}\right) k_{2}D_{2}-2a_{2}\frac{k_{1}}{k_{2}^{2}}D_{1}=0.%
\end{array}%
\right.   \label{elem}
\end{equation}%
Since we wish to analyze the existence of $B_{1}$, $B_{2}$, $D_{1}$, $%
D_{2}\in \left( -\infty ,0\right) $ that solve (\ref{elem}) we couple the
Equations 1 and 3 together%
\begin{equation}
\left( 
\begin{array}{c}
4B_{1}^{2}-\frac{2M_{1}}{k_{1}^{2}} \\ 
4B_{2}^{2}-\frac{2M_{2}}{k_{2}^{2}}%
\end{array}%
\right) =\left( 
\begin{array}{cc}
\frac{2}{k_{1}^{2}}\left( \lambda _{1}+a_{1}\right) k_{1} & -2a_{1}\frac{%
k_{2}}{k_{1}^{2}} \\ 
-2a_{2}\frac{k_{1}}{k_{2}^{2}} & \frac{2}{k_{2}^{2}}\left( \lambda
_{2}+a_{2}\right) k_{2}%
\end{array}%
\right) \left( 
\begin{array}{c}
B_{1} \\ 
B_{2}%
\end{array}%
\right) ,  \label{fixed}
\end{equation}%
and, similarly for the Equations 2 and 4%
\begin{equation}
\left( 
\begin{array}{c}
2B_{1}N \\ 
2B_{2}N%
\end{array}%
\right) =\left( 
\begin{array}{cc}
\frac{2}{k_{1}^{2}}\left( \lambda _{1}+a_{1}\right) k_{1} & -2a_{1}\frac{%
k_{2}}{k_{1}^{2}} \\ 
-2a_{2}\frac{k_{1}}{k_{2}^{2}} & \frac{2}{k_{2}^{2}}\left( \lambda
_{2}+a_{2}\right) k_{2}%
\end{array}%
\right) \left( 
\begin{array}{c}
D_{1} \\ 
D_{2}%
\end{array}%
\right) .  \label{linear}
\end{equation}%
Clearly%
\begin{equation*}
\left\vert 
\begin{array}{cc}
\frac{2}{k_{1}^{2}}\left( \lambda _{1}+a_{1}\right) k_{1} & -2a_{1}\frac{%
k_{2}}{k_{1}^{2}} \\ 
-2a_{2}\frac{k_{1}}{k_{2}^{2}} & \frac{2}{k_{2}^{2}}\left( \lambda
_{2}+a_{2}\right) k_{2}%
\end{array}%
\right\vert =\frac{1}{k_{1}k_{2}}\left( 4\lambda _{1}\lambda _{2}+4\lambda
_{1}a_{2}+4\lambda _{2}a_{1}\right) >0,
\end{equation*}%
and, so the system (\ref{fixed}) can be written equivalently as

\begin{equation}
\left( 
\begin{array}{c}
-B_{1} \\ 
-B_{2}%
\end{array}%
\right) =\frac{1}{2\lambda _{1}\lambda _{2}+2\lambda _{1}a_{2}+2\lambda
_{2}a_{1}}\allowbreak \left( 
\begin{array}{cc}
k_{1}(\lambda _{2}+a_{2}) & \frac{a_{1}k_{2}^{2}}{k_{1}} \\ 
\frac{a_{2}k_{1}^{2}}{k_{2}} & k_{2}(a_{1}+\lambda _{1})%
\end{array}%
\right) \allowbreak \left( 
\begin{array}{c}
\frac{2M_{1}}{k_{1}^{2}}-4B_{1}^{2} \\ 
\frac{2M_{2}}{k_{2}^{2}}-4B_{2}^{2}%
\end{array}%
\right) .  \label{ecs}
\end{equation}%
From the results of \cite{H1,H2} there exist and are unique $B_{1}$, $%
B_{2}\in \left( -\infty ,0\right) $ that solve the system of equations (\ref%
{ecs}). Next, we observe that the system (\ref{linear}) can be written
equivalently as%
\begin{equation*}
\left( 
\begin{array}{c}
D_{1} \\ 
D_{2}%
\end{array}%
\right) =\left( 
\begin{array}{cc}
\frac{\lambda _{2}k_{1}+a_{2}k_{1}}{2\lambda _{1}\lambda _{2}+2\lambda
_{1}a_{2}+2\lambda _{2}a_{1}} & a_{1}\frac{k_{2}^{2}}{2\lambda _{1}\lambda
_{2}k_{1}+2\lambda _{1}a_{2}k_{1}+2\lambda _{2}a_{1}k_{1}} \\ 
a_{2}\frac{k_{1}^{2}}{2\lambda _{1}\lambda _{2}k_{2}+2\lambda
_{1}a_{2}k_{2}+2\lambda _{2}a_{1}k_{2}} & \frac{a_{1}k_{2}^{2}+\lambda
_{1}k_{2}^{2}}{2\lambda _{1}\lambda _{2}k_{2}+2\lambda
_{1}a_{2}k_{2}+2\lambda _{2}a_{1}k_{2}}%
\end{array}%
\right) \allowbreak \left( 
\begin{array}{c}
2B_{1}N \\ 
2B_{2}N%
\end{array}%
\right) ,
\end{equation*}%
from where we can see that there exist $B_{1}$, $B_{2}$, $D_{1}$, $D_{2}\in
\left( -\infty ,0\right) $ that solve (\ref{elem}) and then $\left( 
\underline{u}\left( x\right) ,\underline{v}\left( x\right) \right) $ are
such that the inequalities in (\ref{lam}) hold.

To construct a super-solution it is useful to remember that $\ln 1=0$ and
then a simple calculation shows that 
\begin{equation*}
\text{ }\left( \overline{u}\left( x\right) ,\text{ }\overline{v}\left(
x\right) \right) =\left( 1,1\right) ,
\end{equation*}%
is a super-solution of the problem (\ref{tpe}).

Until now, we constructed the corresponding sub- and super-solutions
employed in the scalar case by \cite{BS}. Clearly, (\ref{tple}) holds and
then in Theorem \ref{mss} it remains to prove that there exists $\left(
u\left( x\right) ,v\left( x\right) \right) \in C^{2}\left( \mathbb{R}%
^{N}\right) \times C^{2}\left( \mathbb{R}^{N}\right) $ with $\underline{u}%
\left( x\right) \leq u\left( x\right) \leq \overline{u}\left( x\right) $ in $%
\mathbb{R}^{N}$ and $\underline{v}\left( x\right) \leq v\left( x\right) \leq 
\overline{v}\left( x\right) $ in $\mathbb{R}^{N}$ satisfying (\ref{tpe}).

To do this, let $B_{k}$ be the ball whose center is the origin of $\mathbb{R}%
^{N}$ and which has radius $k=1,2,...$. We consider the boundary value
problem%
\begin{equation}
\left\{ 
\begin{array}{l}
\Delta u=u[\frac{2}{k_{1}^{2}}\left( f_{1}\left( x\right) +\left( \lambda
_{1}+a_{1}\right) k_{1}\ln u-a_{1}k_{2}\ln v\right) ]\text{, }x\in B_{k}, \\ 
\Delta v=v[\frac{2}{k_{2}^{2}}\left( f_{2}\left( x\right) +\left( \lambda
_{2}+a_{2}\right) k_{2}\ln v-a_{2}k_{1}\ln u\right) ]\text{, }x\in B_{k}, \\ 
u\left( x\right) =\underline{u}_{k}\left( x\right) \text{, }v\left( x\right)
=\underline{v}_{k}\left( x\right) \text{, }x\in \partial B_{k},%
\end{array}%
\right.  \label{ball}
\end{equation}%
where $\underline{u}_{k}=\underline{u}_{\left\vert B_{k}\right. }$ and $%
\underline{v}_{k}=\underline{v}_{\left\vert B_{k}\right. }$. In a similar
way, we define $\overline{u}_{k}=\overline{u}_{\left\vert B_{k}\right. }$
and $\overline{v}_{k}=\overline{v}_{\left\vert B_{k}\right. }$ then $%
\underline{u}_{k}$, $\overline{u}_{k}$, $\underline{v}_{k}$, $\overline{v}%
_{k}\in C^{2}\left( \overline{B}_{k}\right) $.

Observing that%
\begin{eqnarray*}
\underset{x\in \mathbb{R}^{N}}{\inf }\underline{u}\left( x\right) &\leq &%
\underset{x\in \overline{B}_{k}}{\min }\underline{u}_{k}\left( x\right) 
\text{ and }\underset{x\in \mathbb{R}^{N}}{\sup }\overline{u}\left( x\right)
\geq \underset{x\in \overline{B}_{k}}{\max }\overline{u}_{k}\left( x\right) 
\text{,} \\
\underset{x\in \mathbb{R}^{N}}{\inf }\underline{v}\left( x\right) &\leq &%
\underset{x\in \overline{B}_{k}}{\min }\underline{v}_{k}\left( x\right) 
\text{ and }\underset{x\in \mathbb{R}^{N}}{\sup }\overline{v}\left( x\right)
\geq \underset{x\in \overline{B}_{k}}{\max }\overline{v}_{k}\left( x\right) 
\text{,}
\end{eqnarray*}%
a result of Reis Gaete \cite{MG} (see also the pioneering papers of Kawano 
\cite{NK} and Lee, Shivaji and Ye \cite{LSY}), proves the existence of a
solution $\left( u_{k},v_{k}\right) \in \left[ C^{2}\left( B_{k}\right) \cap
C\left( \overline{B}_{k}\right) \right] ^{2}$ satisfying the system (\ref%
{ball}). The functions $\left( u_{k},v_{k}\right) $ also satisfy%
\begin{eqnarray*}
\underline{u}_{k}\left( x\right) &\leq &u_{k}\left( x\right) \leq \overline{u%
}_{k}\left( x\right) \text{, }x\in \overline{B}_{k}, \\
\underline{v}_{k}\left( x\right) &\leq &v_{k}\left( x\right) \leq \overline{v%
}_{k}\left( x\right) \text{, }x\in \overline{B}_{k}.
\end{eqnarray*}%
By a standard regularity argument based on Schauder estimates, see Tolksdorf 
\cite[17, proposition 3.7, p. 806]{T} and Reis Gaete \cite{MG} for details,
we can see that for all integers $k\geq n+1$ there are $\alpha _{1},\alpha
_{2}\in \left( 0,1\right) $ and positive constants $C_{1}$, $C_{2},$
independent of $k$ such that 
\begin{equation}
\left\{ 
\begin{array}{c}
u_{k}\in C^{2,\alpha _{1}}\left( \overline{B}_{n}\right) \text{ and }%
\left\vert u_{k}\right\vert _{C^{2,\alpha _{1}}\left( \overline{B}%
_{n}\right) }<C_{1}, \\ 
v_{k}\in C^{2,\alpha _{2}}\left( \overline{B}_{n}\right) \text{ and }%
\left\vert v_{k}\right\vert _{C^{2,\alpha _{2}}\left( \overline{B}%
_{n}\right) }<C_{2},%
\end{array}%
\right.  \label{b1}
\end{equation}%
where $\left\vert \circ \right\vert _{C^{2,\circ }}$ is the usual norm of
the space $C^{2,\circ }\left( \overline{B}_{n}\right) $. Moreover, there
exist constants: $C_{3}$ independent of $u_{k}$, $C_{4}$ independent of $%
v_{k}$ and such that%
\begin{equation}
\left\{ 
\begin{array}{c}
\underset{x\in \overline{B}_{n}}{\max }\left\vert \nabla u_{k}\left(
x\right) \right\vert \leq C_{3}\underset{x\in \overline{B}_{k}}{\max }%
\left\vert u_{k}\left( x\right) \right\vert , \\ 
\underset{x\in \overline{B}_{n}}{\max }\left\vert \nabla v_{k}\left(
x\right) \right\vert \leq C_{4}\underset{x\in \overline{B}_{k}}{\max }%
\left\vert v_{k}\left( x\right) \right\vert .%
\end{array}%
\right.  \label{b2}
\end{equation}%
The information from (\ref{b1}) and (\ref{b2}) implies that \{$\left( \nabla
u_{k},\nabla v_{k}\right) $\}$_{k}$ as well as \{$\left( u_{k},v_{k}\right) $%
\}$_{k}$ are uniformly bounded on $\overline{B}_{n}$. We wish to show that
this sequence \{$\left( u_{k},v_{k}\right) $\}$_{k}$ contains a subsequence
converging to a desired entire solution of (\ref{tpe}). Next, we concentrate
our attention to the sequence \{$u_{k}$\}$_{k}$. Using the compactness of
the embedding $C^{2,\alpha _{1}}\left( \overline{B}_{n}\right)
\hookrightarrow C^{2}\left( \overline{B}_{n}\right) $, enables us to define
the subsequence%
\begin{equation*}
u_{n}^{k}:=u_{k\left\vert B_{n}\right. },\text{ for all }k\geq n+1\text{.}
\end{equation*}%
Then for $n=1,2,3,...$ there exist a subsequence $\{u_{n}^{k_{nj}}\}_{k\geq
n+1,j\geq 1}$ of $\{u_{n}^{k}\}_{k\geq n+1}$ and a function $u_{n}$ such that%
\begin{equation}
u_{n}^{k_{nj}}\rightarrow u_{n}\text{,}  \label{sub}
\end{equation}%
uniformly in the $C^{2}\left( \overline{B}_{n}\right) $ norm. More exactly,
we get through a well-known diagonal process that%
\begin{eqnarray*}
\mathbf{u}_{1}^{k_{11}}\text{, }u_{1}^{k_{12}}\text{, }u_{1}^{k_{13}}\text{, 
}... &\longrightarrow &u_{1}\text{ in }C^{2}\left( \overline{B}_{1}\right) ,
\\
u_{2}^{k_{21}}\text{, }\mathbf{u}_{2}^{k_{22}}\text{, }u_{2}^{k_{23}}\text{, 
}... &\longrightarrow &u_{2}\text{ in }C^{2}\left( \overline{B}_{2}\right) ,
\\
u_{3}^{k_{31}}\text{, }u_{3}^{k_{32}}\text{, }\mathbf{u}_{3}^{k_{33}}\text{, 
}... &\longrightarrow &u_{3}\text{ in }C^{2}\left( \overline{B}_{3}\right) ,
\\
&&...
\end{eqnarray*}%
Since $\mathbb{R}^{N}=\underset{n=1}{\overset{\infty }{\cup }}B_{n}$, we can
define the function $u:\mathbb{R}^{N}\rightarrow \left[ 0,\infty \right) $
such that 
\begin{equation*}
u\left( x\right) =\lim_{n\rightarrow \infty }u_{n}\left( x\right) .
\end{equation*}%
Let us give the construction of the function $u$ for the problem (\ref{tpe}%
). This is obtained by considering the sequence ($u_{d}^{k_{dd}}$)$_{d\geq
1} $ and the sequence $(u_{n}^{k_{nd}})_{k\geq n+1}$, restricted to the ball 
$B_{n}$, which are such that%
\begin{equation*}
u_{n}^{k_{nd}}\overset{d\rightarrow \infty }{\rightarrow }u_{n}:=u\left(
x\right) \text{ for all }x\in B_{n}\text{,}
\end{equation*}%
and then, for $d\rightarrow \infty $ we obtain%
\begin{equation*}
u_{d}^{k_{dd}}\overset{d\rightarrow \infty }{\rightarrow }u\left( x\right) 
\text{ in }C^{2}\left( \mathbb{R}^{N}\right) \text{,}
\end{equation*}%
according with the diagonal process. Furthermore, since%
\begin{equation*}
\underline{u}\left( x\right) \leq u_{d}^{k_{dd}}\leq \overline{u}\left(
x\right) \text{, for }x\in \mathbb{R}^{N},
\end{equation*}%
and for each $d=1,2,3,...$ the following relation is valid 
\begin{equation*}
\underline{u}\left( x\right) \leq u\left( x\right) \leq \overline{u}\left(
x\right) \text{, for }x\in \mathbb{R}^{N}.
\end{equation*}%
We employ the same iteration scheme to construct the function $v:\mathbb{R}%
^{N}\rightarrow \left[ 0,\infty \right) $ such that 
\begin{equation*}
v\left( x\right) =\lim_{n\rightarrow \infty }v_{n}\left( x\right) .
\end{equation*}%
From the regularity theory the solution $\left( u,v\right) $ belongs to $%
C^{2}\left( \mathbb{R}^{N}\right) \times C^{2}\left( \mathbb{R}^{N}\right) $
and satisfies (\ref{tpe}). This completes the proof of Theorem \ref{mss}.

\paragraph{Proof of Theorem \protect\ref{msg}}

As easily verified, the existence of solutions is proved by Lemma \ref{ech}
and Theorem \ref{mss}. Then it remains to prove (\ref{ineq}).

A recapitulation of the changes of variables say that 
\begin{equation}
z_{1}\left( x\right) =-k_{1}\ln u\left( x\right) \text{ and }z_{2}\left(
x\right) =-k_{2}\ln v\left( x\right) \text{,}  \label{neg}
\end{equation}%
is a solution for (\ref{tpes}). Observing that%
\begin{equation*}
\underline{u}\left( x\right) =e^{B_{1}\left\vert x\right\vert
^{2}+D_{1}}\leq u\left( x\right) \leq \text{ }\overline{u}\left( x\right) =1%
\text{, }x\in \mathbb{R}^{N},
\end{equation*}%
it follows that%
\begin{equation*}
B_{1}\left\vert x\right\vert ^{2}+D_{1}\leq \ln u\left( x\right) \leq \ln 1,
\end{equation*}%
and, then%
\begin{equation*}
0\leq -k_{1}\ln u\left( x\right) \leq -k_{1}B_{1}\left\vert x\right\vert
^{2}-k_{1}D_{1},
\end{equation*}%
or equivalently%
\begin{equation*}
0\leq z_{1}\left( x\right) \leq K_{1}\left( \left\vert x\right\vert
^{2}+1\right) \text{, for }x\in \mathbb{R}^{N}\text{ and }K_{1}=\max
\{-k_{1}B_{1},-k_{1}D_{1}\}.
\end{equation*}%
In the same way%
\begin{equation*}
0\leq z_{2}\left( x\right) \leq K_{2}\left( \left\vert x\right\vert
^{2}+1\right) \text{, for }x\in \mathbb{R}^{N}\text{ and }K_{2}=\max
\{-k_{2}B_{2},-k_{2}D_{2}\},
\end{equation*}%
and the proof is completed.

By the same arguments as in \cite[Theorem 3, p. 278]{AL} the solution $%
\left( z_{1}\left( x\right) ,z_{2}\left( x\right) \right) $ is convex. Since 
$\left( z_{1}\left( x\right) ,z_{2}\left( x\right) \right) $ verifies (\ref%
{ineq}) the inequality (\ref{sass}) follows from \cite[Lemma 1, p. 24]{E}
(see also the arguments in \cite[Theorem 1, p. 236]{EP}). The uniqueness of
such a solution follows from the result of \cite{A,HE,KAT} (see also the
former papers of \cite{ISHI1,ISHI2}), since for our system their comparision
results can also be set in $\mathbb{R}^{N}$, instead of a domain $\Omega
\subset \mathbb{R}^{N}$. The proof is completed.

\section{Context and the Diffusion Model \label{ppp}}

Let us present the setting. Consider $W$ a $N-$dimensional Brownian motion
on a filtered probability space 
\begin{equation}
(\Omega ,\{\mathcal{F}_{t}\}_{0\leq t\leq T},\mathcal{F},P),  \label{000}
\end{equation}%
where $\{\mathcal{F}_{t}\}_{0\leq t\leq T}$ is a completed filtration, and $%
T=\infty $ (we deal with the infinite horizon case). We allow for regime
switching in our model; regime switching refers to the situation when the
characteristics of the state process are affected by several regimes (e.g.
in finance bull and bear market with higher volatility in the bear market).
The regime switching is captured by a continuous time homogeneous Markov
chain $\epsilon (t)$ adapted to $\mathcal{F}_{t}$ with two regimes good and
bad, i.e., for every 
\begin{equation*}
t\in \lbrack 0,\infty )\text{ and }\epsilon (t)\in \{ {1}, {2 } \}.
\end{equation*}%
In a specific application, $\epsilon (t)=1$ could represent a regime of
economic growth while $\epsilon (t)=2$ could represent a regime of economic
recession. In another application, $\epsilon (t)=1$ could represent a regime
in which consumer demand is high while $\epsilon (t)=2$ could represent a
regime in which consumer demand is low.

The Markov chain's rate matrix is 
\begin{equation}
A=\left( 
\begin{array}{cc}
-a_{1} & a_{1} \\ 
a_{2} & -a_{2}%
\end{array}%
\right) ,  \label{4}
\end{equation}%
for some $a_{1}>0,$ $a_{2}>0$. Diagonal elements $A_{ii}$ are defined such
that%
\begin{equation}
A_{ii}=-\underset{j\neq i}{\Sigma }A_{ij},  \label{5}
\end{equation}%
where 
\begin{equation*}
A_{11}=-a_{1},A_{12}=a_{1},A_{21}=a_{2},A_{22}=-a_{2}.
\end{equation*}%
In this case, if $p_{t}=\mathbb{E}[\epsilon (t)]\in \mathbb{R}^{2},$ then 
\begin{equation}
\frac{d\epsilon (t)}{dt}=A\epsilon (t).  \label{6}
\end{equation}%
Moreover 
\begin{equation}
\epsilon (t)=\epsilon (0)+\int_{0}^{t}A\epsilon (u)\,du+M_{t},  \label{mj}
\end{equation}%
where ${M(t)}$ is a martingale with respect to $\mathcal{F}_{t}.$

Let us consider a Markov modulated controlled diffusion with controls in
feed-back form 
\begin{equation}
dX^{i}(t)=c^{i}_{\epsilon (t)}(X(t))dt+k_{\epsilon (t)} dW^{i}(t),
i=1,\ldots N,  \label{sde}
\end{equation}%
for some constants $k_{1}>0,$ $k_{2}>0,$ and $X(0)=x\in \mathbb{R}^{N}$. \
Here, at every time $t$, the control $c_{\epsilon (t)} $ (for instance the
demand of certain items), and the volatility $k_{\epsilon (t)}$ depend on
the regime $\epsilon (t)$. We allow the demand to take on negative values,
which represent items return (due to spoilage). We consider the class of
admissible controls, $\mathcal{A} $, which are the feedback controls for
which the SDE \eqref{sde} has a unique strong solution.

The infinitesimal generator $L$ of diffusion $X$ is second order
differential operator defined by 
\begin{equation}
L^{c} v(x, 1 )= \frac{1}{2} k_{ 1 } \Delta v (x, 1 )+ c_{ 1 } \nabla v (x, 1
) +A_{ 1 1} v(x, 1 ) + A_{ 1 2} v(x,2 ),  \label{30}
\end{equation}

\begin{equation}
L^{c}v(x,2)=\frac{1}{2}k_{2}\Delta v(x,2)+c_{2}\nabla
v(x,2)+A_{22}v(x,2)+A_{21}v(x,1),  \label{31}
\end{equation}%
(see \cite{M} for more on this). Following this we can state Itô's formula
for Markov modulated diffusion

\begin{equation}
dv(X(t),\epsilon (t))=L^{c}v(X(t),\epsilon (t))dt+k_{\epsilon (t)}\nabla
v(X(t),\epsilon (t))dW(t).  \label{32}
\end{equation}%
Next, for each $c\in \mathcal{A}$ the cost functional is defined by 
\begin{equation}
J(x,c,i)=E[\int_{0}^{\infty }e^{-\lambda _{\epsilon (t)}t}[f_{\epsilon
(t)}(X(t))+\frac{1}{2}|c|_{\epsilon (t)}^{2}(X(t))]\,dt|\epsilon (0)=i].
\label{7}
\end{equation}%
Our objective is to minimize the functional $J$, i.e. determine the value
function 
\begin{equation}
z_{i}(x)=\inf {J(x,c,i)},  \label{8}
\end{equation}%
and to find the optimal control. The infimum is taken over all admissible
controls $c\in \mathcal{A}.$ Notice that the discount rate depends on the
regime; for more on this modelling approach see \cite{PZ}.

In order, to obtain the HJB equation, we apply the
martingale/supermartingale principle; search for a function $u(x,i)$ such
that the stochastic process $M^{c}(t)$ defined below 
\begin{equation}
M^{c}(t)=e^{-\lambda _{\epsilon (t)}}u(X(t),\epsilon
(t))-\int_{0}^{t}e^{-\lambda _{\epsilon (u)}u}[f_{\epsilon (t)}(X(u))+\frac{1%
}{2}|c|_{\epsilon (u)}^{2}(X(u))]\,du,  \label{99}
\end{equation}%
is supermartingale and martingale for the optimal control. If this is
achieved together with the following transversality condition 
\begin{equation}
\lim_{t\rightarrow \infty }E[e^{-\lambda _{\epsilon (t)}t}u(X(t),\epsilon
(t))]=0,  \label{tr}
\end{equation}%
and some estimates on the value function yield that 
\begin{equation}
z_{i}(x)=-u(x,i)=\inf_{c\in \mathcal{A}}{J(x,c,i)}.  \label{10}
\end{equation}%
The proof of this statement is done in the Verification subsection.

The supermartingale/martingale requirement leads to the following HJB
equation 
\begin{equation}
\frac{k_{i}}{2}\Delta u(x,i)+\sup_{c\in \mathcal{A}}[\nabla u(x,i)c-\frac{%
|c|^{2}}{2}]=f_{i}\left( x\right) +(\lambda _{i}+a_{i})u(x,i)-a_{i}u(x,j)%
\mathbf{,}  \label{HJB1}
\end{equation}%
for $\,i$, $j\in \{1,2\}$. First order condition yields the candidate
optimal control 
\begin{equation}
\hat{c}_{i}(x)=\nabla u(x,i)=-\nabla z_{i}(x),  \label{12}
\end{equation}%
and this leads to the system 
\begin{equation}
\frac{k_{i}}{2}\Delta u(x,i)+\frac{|\nabla u(x,i)|^{2}}{2}=f_{i}(x)+(\lambda
_{i}+a_{i})u(x,i)-a_{i}u(x,j),\,  \label{HJB11}
\end{equation}%
for $\,i$, $j\in \{{1},{2}\}$. Alternatively this system can be written in
terms of $z_{i}(x)$, ($i=1,2$) to get (\ref{tpes}).

\section{Verification \label{avr}}

In this section we establish the optimality of control 
\begin{equation}
\hat{c}_{i}(x)=\nabla u(x,i)=-\nabla z_{i}(x)\text{.}  \label{012}
\end{equation}%
Its associated Markov modulated diffusion is 
\begin{equation}
dX^{i}(t)=\hat{c}_{\epsilon (t)}^{i}(X(t))dt+k_{\epsilon
(t)}dW^{i}(t),i=1,\ldots N.  \label{sde0}
\end{equation}

The verification theorem proceeds with the following steps:

\textbf{First Step:} Girsanov theorem for Markov-modulated processes (Lemma
1 page 286 in \cite{YZY}) together with \eqref{sass} yield a weak solution
for SDE (\ref{sde0}).

\textbf{Second Step:} Let $X(t)$ be the solution of (\ref{sde0}). In light
of (\ref{sass}) one can get using exercise $7.5$ of \cite{Oks} that 
\begin{equation}
E\left\vert X(t)\right\vert ^{2}\leq C_{1}e^{C_{2}t},  \label{est}
\end{equation}%
for some positive constants $C_{1},{C_{2}.}$

\textbf{Third Step:} The set of acceptable controls that we consider is
encompassing of controls $c$ for which 
\begin{equation}
J(x,c,i)=E[\int_{0}^{\infty }e^{-\lambda _{\epsilon (t)}t}[f_{\epsilon
(t)}(X(t))+\frac{1}{2}|c|_{\epsilon (t)}^{2}(X(t))]\,dt|\epsilon
(0)=i]<\infty ,  \label{70}
\end{equation}%
and the following transversality condition%
\begin{equation*}
\lim_{t\rightarrow \infty }Ee^{-\lambda _{\epsilon (t)}}\left\vert
X(t)\right\vert ^{2}=0,
\end{equation*}%
is met. Because of (\ref{sass}), estimates (\ref{ineq}), (\ref{est}), the
candidate optimal control $\hat{c}$ of (\ref{012}) verifies that $%
J(x,c,i)<\infty $, for $\lambda _{1},\lambda _{2}$ large enough. Moreover,
there exist $\lambda _{1}>0$ and $\lambda _{2}>0$ large enough such that the
transversality condition (\ref{tr}) is met because of (\ref{ineq}) and (\ref%
{est}). Also the control $c=0$, is an acceptable control.

In light of the quadratic estimate on the value function (see (\ref{ineq})
in theorem 2.1), the transversality condition implies that 
\begin{equation}  \label{trans1}
\lim_{t\rightarrow \infty }Ee^{-\lambda_{\epsilon (t)}} u ( X(t), \epsilon
(t)) =0.
\end{equation}

\textbf{Fourth Step:} Recall that 
\begin{equation}
M^{c}(t)=e^{-\lambda _{\epsilon (t)}}u(X(t),\epsilon
(t))-\int_{0}^{t}e^{-\lambda _{\epsilon (u)}u}[f_{\epsilon (t)}(X(u))+\frac{1%
}{2}|c|_{\epsilon (u)}^{2}(X(u))]\,du.  \label{9}
\end{equation}%
Therefore, the Itô's Lemma yields for the optimal control candidate, $\hat{c}
$%
\begin{equation*}
dM^{c}\left( t\right) =e^{-\lambda _{\epsilon (t)}}k_{\epsilon (t)}\nabla
u(X(t),\epsilon (t))dW(t).
\end{equation*}%
Consequently $M^{\hat{c}}(t)$ is a local martingale. Moreover, for $\lambda
_{1},\lambda _{2}$ large enough, in light of (\ref{sass}), and (\ref{est}), 
\begin{equation*}
E\int_{0}^{t}e^{-2\lambda _{\epsilon (s)}}k_{\epsilon (s)}^{2}\left\vert
\nabla u(X(s),\epsilon (s))\right\vert ^{2}ds\leq \overline{C},
\end{equation*}%
for some positive constants $\overline{C}$. This in turn makes $M^{\hat{c}%
}(t)$ a (true) martingale.

\textbf{Fifth Step: }This step establishes the optimality of $\hat{c}$ of (%
\ref{012}). The HJB equation (\ref{HJB1}) is equivalent to 
\begin{equation*}
\sup_{c}L^{c}u(x,i)=0,\quad L^{\hat{c}}u(x,i)=0\text{, }i=1,2.
\end{equation*}

The martingale/supermartingale principle yields%
\begin{equation*}
Ee^{-\lambda _{\epsilon (t)}}u(X(t),\epsilon (t))-E\int_{0}^{t}e^{-\lambda
_{\epsilon (u)}u}[f_{\epsilon (t)}(X(u))+\frac{1}{2}|\hat{c}|_{\epsilon
(u)}^{2}(X(u))]\,du=u(x,\epsilon (0)),
\end{equation*}%
and%
\begin{equation*}
Ee^{-\lambda _{\epsilon (t)}}u(X(t),\epsilon (t))-E\int_{0}^{t}e^{-\lambda
_{\epsilon (u)}u}[f_{\epsilon (t)}(X(u))+\frac{1}{2}|{c}|_{\epsilon
(u)}^{2}(X(u))]\,du\leq u(x,\epsilon (0)).
\end{equation*}%
By passing $t\rightarrow \infty $ and using transversality condition (\ref%
{trans1}) we get the optimality of $\hat{c}$.

\section{Special Case \label{sc}}

In the following we manage to obtain a simple closed form solution for our
system given a special discount $\lambda _{1}$, $\lambda _{2}$ and the loss
functions of the type $f_{1}(x)=f_{2}\left( x\right) =\left\vert
x\right\vert ^{2}$. That is, assume 
\begin{equation*}
\begin{array}{c}
\lambda _{1}=-a_{1}+Nk_{1}+\frac{1}{8}a_{1}\left( \sqrt{N^{2}k_{1}^{2}+8}%
-Nk_{1}\right) \left( \sqrt{N^{2}k_{2}^{2}+8}-Nk_{2}\right) \\ 
+\frac{1}{4}Na_{1}k_{1}\left( \sqrt{N^{2}k_{2}^{2}+8}-Nk_{2}\right) , \\ 
\lambda _{2}=-a_{2}+Nk_{2}+\frac{1}{8}a_{2}\left( \sqrt{N^{2}k_{1}^{2}+8}%
-Nk_{1}\right) \left( \sqrt{N^{2}k_{2}^{2}+8}-Nk_{2}\right) \\ 
+\frac{1}{4}Na_{2}k_{2}\left( \sqrt{N^{2}k_{1}^{2}+8}-Nk_{1}\right) ,%
\end{array}%
\end{equation*}%
are such that $\lambda _{1}>0,\lambda _{2}>0$. Then, one solution for the
problem (\ref{tpe}) is%
\begin{eqnarray*}
u\left( \left\vert x\right\vert \right) &=&e^{m_{1}\left( \left\vert
x\right\vert ^{2}+1\right) }\text{, }m_{1}=-\frac{1}{4k_{1}}\left( \sqrt{%
N^{2}k_{1}^{2}+8}-Nk_{1}\right) \text{,} \\
\text{ }v\left( \left\vert x\right\vert \right) &=&e^{m_{2}\left( \left\vert
x\right\vert ^{2}+1\right) }\text{, }m_{2}=-\frac{1}{4k_{2}}\left( \sqrt{%
N^{2}k_{2}^{2}+8}-Nk_{2}\right) .
\end{eqnarray*}%
Let us point out that (\ref{neg}) implies%
\begin{equation}
z_{1}\left( x\right) =-k_{1}m_{1}\left( \left\vert x\right\vert
^{2}+1\right) >0\text{ and }z_{2}\left( x\right) =-k_{2}m_{2}\left(
\left\vert x\right\vert ^{2}+1\right) >0\text{ for all }x\in \mathbb{R}^{N}%
\text{,}  \label{poz1}
\end{equation}%
i.e. $\left( z_{1}\left( x\right) ,z_{2}\left( x\right) \right) $ is the
positive solution obtained with the above procedure. For the stochastic
control problem we choose the positive solution, i.e., the one given in (\ref%
{poz1}).

We also point that in the special case $f_{i}\left( x\right)
=B_{i}\left\vert x\right\vert ^{2}+C_{i}\left\vert x\right\vert +D_{i}$,
with $A_{i}$, $B_{i}$, $C_{i}\in \mathbb{R}$ suitable chosen, many others
exact solutions of the form 
\begin{equation*}
\left( z_{1}\left( \left\vert x\right\vert \right) ,z_{2}\left( x\right)
\right) =\left( A_{1}\left\vert x\right\vert ^{2}+A_{2}\left\vert
x\right\vert +A_{3},A_{4}\left\vert x\right\vert ^{2}+A_{5}\left\vert
x\right\vert +A_{6}\right) \text{, }A_{i}\in \mathbb{R}_{+}\text{,}
\end{equation*}%
can be constructed using the computation technique. Since the solution form
is more complicated we omit to give his closed form. Moreover, by our
Theorem \ref{msg} such a solution is unique.

\textbf{Acknowledgments.} Traian A. Pirvu acknowledges that this work was
supported by NSERC grant 5-36700.

\bigskip

\end{document}